\documentclass[preprint,12pt]{elsarticle}

\journal{Discrete Applied Mathematics}
\usepackage{amsthm}
\usepackage{enumerate}
\usepackage{amssymb}
\usepackage{amsmath}
\usepackage{amsfonts}

\newtheorem{thrm}{Theorem}[section]
\newcommand{\thm}{\begin{thrm}}
\newcommand{\xthm}{\end{thrm}}
\newtheorem{corl}[thrm]{Corollary}
\newcommand{\cor}{\begin{corl}}
\newcommand{\xcor}{\end{corl}}
\newtheorem{lemm}[thrm]{Lemma}
\newcommand{\lem}{\begin{lemm}}
\newcommand{\xlem}{\end{lemm}}
\newtheorem{claim}[thrm]{Claim}
\newcommand{\clm}{\begin{claim}}
\newcommand{\xclm}{\end{claim}}
\newtheorem{propos}[thrm]{Proposition}
\newcommand{\prop}{\begin{propos}}
\newcommand{\xprop}{\end{propos}}
\newtheorem{conjec}[thrm]{Conjecture}
\newcommand{\conj}{\begin{conjec}}
\newcommand{\xconj}{\end{conjec}}
\newcommand{\prf}{\begin{proof}}
\newcommand{\xprf}{\end{proof}}

\newcommand{\sm}{\setminus}

\newcommand{\dta} {\delta}

\newcommand{\lam} {\lambda}

\newcommand{\sig} {\sigma}

\newcommand{\Gam}{\Gamma}

\newcommand{\arr}{\begin{array}}
\newcommand{\xarr}{\end{array}}
\newcommand{\tabl}{\begin{tabular}}
\newcommand{\xtabl}{\end{tabular}}

\newcommand{\f}{\frac}

\newcommand{\eqns}{\begin{eqnarray*}}
\newcommand{\xeqns}{\end{eqnarray*}}
\newcommand{\enum}{\begin{enumerate}}
\newcommand{\xenum}{\end{enumerate}}
\newcommand{\itmz}{\begin{itemize}}
\newcommand{\xitmz}{\end{itemize}}
\newcommand{\case}{\begin{cases}}
\newcommand{\xcase}{\end{cases}}

\newcommand{\inv}{^{-1}}

\newcommand{\bmat}{\begin{bmatrix}}
\newcommand{\xbmat}{\end{bmatrix}}
\newcommand{\pmat}{\begin{pmatrix}}
\newcommand{\xpmat}{\end{pmatrix}}

\newcommand{\prob}{\begin{problem}}
\newcommand{\xprob}{\end{problem}}
\newcommand{\soln}{\begin{solution}}
\newcommand{\xsoln}{\end{solution}}

\begin{document}
\date{}

\begin{frontmatter}
\title{On Coupon Colorings of Graphs}
\author[ucsd]{Bob Chen}\ead{b2chen@math.ucsd.edu}
\author[kias]{Jeong Han Kim}\ead{hmkkim@gmail.com}
\author[ucsd]{Michael Tait \corref{cor1}}\ead{mtait@math.ucsd.edu}\cortext[cor1]{Corresponding author}
\author[ucsd]{Jacques Verstraete\fnref{fn1}}\ead{jverstra@math.ucsd.edu}
\fntext[fn1]{Research supported by NSF Grant DMS 1101489.}
\address[ucsd]{University of California, San Diego}
\address[kias]{Korean Institute for Advanced Study}

\begin{abstract}
Let $G$ be a graph with no isolated vertices. A {\em $k$-coupon coloring} of $G$ is an assignment of colors from $[k] := \{1,2,\dots,k\}$ to the vertices of $G$ such that
the neighborhood of every vertex of $G$ contains vertices of all colors from $[k]$. The maximum $k$ for which a $k$-coupon coloring exists is called
the {\em coupon coloring number} of $G$, and is denoted $\chi_{c}(G)$. In this paper, we prove that every $d$-regular graph $G$
has $\chi_{c}(G) \geq (1 - o(1))d/\log d$ as $d \rightarrow \infty$, and the proportion of $d$-regular graphs $G$ for which
$\chi_c(G) \leq (1 + o(1))d/\log d$ tends to $1$ as $|V(G)| \rightarrow \infty$.
\end{abstract}


\end{frontmatter}
\section{Introduction}

Let $G$ be a graph with no isolated vertices. A {\em $k$-coupon coloring} of $G$ is an assignment of colors from $[k] := \{1,2,\dots,k\}$ to the vertices of $G$ such that the neighborhood of every vertex of $G$ contains vertices of all colors from $[k]$. The maximum $k$ for which a $k$-coupon coloring exists is called the {\em coupon coloring number} of $G$, and is denoted $\chi_{c}(G)$. The coupon coloring number of any graph $G$ with no isolated vertices is well-defined, since we may assign every vertex the same color.

\medskip

The motivation for the term coupon coloring is that we may imagine the colors as coupons of
 different types, and then the requirement of coupon coloring is that every vertex collects from its neighbors coupons of all different types.
 If we imagine that users $v_1,v_2,\dots,v_n$ are each assigned a bit from a $k$-bit message, and that every user has contact with a set of other users,
 then every user can reconstruct from her contacts the entire message if and only if the graph of contacts has a $k$-coupon coloring.
 The task given a graph of contacts is to determine the coupon coloring number of the graph, to maximize the length of the message that can be transmitted.
 It is possible to give examples of graphs of very large minimum degree whose coupon coloring number is 1; for this reason, we consider $d$-regular graphs in this paper.
 The coupon coloring problem for hypercubes is closely related to problems in coding theory (see \"{O}sterg\aa rd~\cite{O}), and for $k = 2$ the coupon coloring problem is also equivalent to the well-studied Property B of hypergraphs~\cite{E}: if we form from a graph $G$ the {\em neighborhood hypergraph} $H = \{\Gamma(v) : v \in V(G)\}$ with vertex set $V(G)$,
then a 2-coupon coloring of $G$ exists if and only if $H$ has Property B -- namely, a 2-coloring of the vertices of $H$ such that
no edge of $H$ is monochromatic. More generally, if $H$ is a hypergraph and $H$ admits a $k$-coloring such that every edge contains all
$k$ colors, then $H$ is said to have a {\em panchromatic $k$-coloring}, which corresponds to a $k$-coupon coloring. For fixed $k$,
local conditions on $H$ were given by Kostochka and Woodall~\cite{KW} for a hypergraph to have a panchromatic $k$-coloring. In particular
it is shown that if every edge in $H$ has size at least $k$, and for every set $E$ of edges of $H$ one has $|\bigcup E| \geq (k - 1)|E| - k + 3$, then  $H$ has a panchromatic $k$-coloring. In this paper we study the extremal values of $\chi_c(G)$ when $G$ is a $d$-regular graph as $d$ becomes large.

In addition, results on coupon colorings have concrete applications in network science. One application is to large multi-robot networks (cf \cite{BCM}). One may imagine a network large enough that robots must act based on local information. A 
graph can be constructed with robots in the network as nodes and an edge between nodes if the corresponding robots are able 
to communicate with each other. An example described in \cite{AELTW} is as follows: a group of robots is deployed to 
monitor an environment. Each robot must monitor many different statistics (e.g. temperature, humidity, etc.), but due to 
power limitations is only equipped with a single sensor (thermometer, barometer, etc.). Thus, in order to obtain the 
remaining data, each robot must communicate with its neighbors.  A similar example arises in allocating resources to a 
network \cite{AELTW}. If each vertex of a graph may only use resources available at the vertex or its neighbors, and if 
some resource (e.g. a printer) must be available to every node in the network, then copies of that resource must be 
allocated to a dominating set of the network. If every node in the network can accomodate one resource, then finding the 
coupon coloring number of the network is equivalent to finding the maximum number of resources that can be made available 
to every node in the network.

\subsection{Coupon coloring regular graphs}

For every positive integer $d$, there exists a graph of minimum degree at least $d$ with coupon coloring number equal to 1: for instance
this is easily verified for the bipartite graph of minimum degree $d$ which represents the incidence graph of a complete $d$-uniform hypergraph on more than $d(d - 1)$ vertices, since this
hypergraph does not have Property B. Therefore $\chi_c(G)$ cannot be controlled by the minimum degree of a graph. The main theorem of this paper
determines the asymptotic value of the coupon coloring number of almost all $d$-regular graphs,
and also shows that every $d$-regular graph has coupon coloring number at least about $d/\log d$. Let $G_{n,d}$ denote
the probability space of $d$-regular $n$-vertex graphs, with the uniform probability measure.

\thm \label{main}
For every $\delta > 0$, there exists $d_0(\delta)$ such that if $d \geq d_0(\epsilon)$, then every $d$-regular graph $G$ has
\[ \chi_{c}(G) \geq (1 - \delta)\frac{d}{\log d}.\]
For every $\epsilon > 0$, there exists $d_1(\epsilon)$ such that if $d \geq d_1(\epsilon)$, then as $n \rightarrow \infty$, almost every $d$-regular $n$-vertex graph $G$ has
\[ \chi_c(G) \leq (1 + \epsilon)\frac{d}{\log d}.\]
\xthm

\medskip

In particular, this theorem says $\chi_c(G) \sim d/(\log d)$ as $d \rightarrow \infty$ for almost every $d$-regular graph $G$.
Theorem \ref{main} is proved in Section 4. An explicit example of a regular
graph with small coupon coloring number is the {\em Paley graph} $G_q$ over $\mathbb F_q$ when $q \equiv 1$ mod 4 is a prime power. This graph
is a $d = (q - 1)/2$-regular graph with $q$ vertices. The minimum size of a dominating set in $G_q$ is known to be at least $(\frac{1}{2} - o(1))\log q$, using standard character sum estimates (see G\'{a}cs and Sz\H{o}nyi~\cite{GS}), and since a coupon coloring requires every color class to be a
dominating set, $\chi_c(G_q) \leq (1 + o(1))2q/\log q \leq (4 + o(1))d/\log d$, which is within a factor four of the bound in Theorem \ref{main}.
It would be interesting to provide explicit constructions of $d$-regular graphs with coupon coloring number $(1 + o(1))d/\log d$, and in particular
to determine $\chi_c(G_q)$ when $G_q$ denotes the Paley graph over $\mathbb F_q$.

\bigskip

\subsection{Coloring cubes}

We next consider graphs whose coupon coloring number is as large as possible. It is not hard to give examples of $d$-regular graphs $G$ with $\chi_c(G) = d$, for instance, the complete bipartite
graph with $d$ vertices in each part has coupon coloring number equal to $d$. An {\em injective $k$-coloring} of a graph $G$ is an assignment of colors from $[k]$ to the vertices of $G$ such that no two vertices joined by a path of length two in $G$ have the same color.
The minimum $k$ for which such a coloring exists is called the {\em injective coloring number} of $G$, denoted $\chi_{i}(G)$, and is well-defined since we may assign all vertices of $G$ different colors.
One observes
\[
\chi_c(G) \leq \delta(G) \leq \Delta(G) \leq \chi_{i}(G),
\]
 where $\delta(G)$ and $\Delta(G)$ denote the minimum and maximum degree of $G$. If $G$ is $d$-regular and $\chi_{i}(G) = d$ then also
 $\chi_c(G) = d$, and this is the largest possible value of $\chi_c(G)$.

 \bigskip

 The injective coloring number
 has been studied for numerous classes of graphs, especially planar graphs~\cite{Bu,Bu2,C,C2,Luzar}.
 The structure of $d$-regular graphs with $\chi_{i}(G) = d = \chi_c(G)$ can be described as follows:
 they may be constructed from a complete graph of order $d$ by replacing each vertex with a set of $n/d$ vertices, placing a perfect matching on these sets of $n/d$ vertices
 and placing a perfect matching between any pair of distinct sets of those $n/d$ vertices. In particular, in a $d$-coupon coloring all color classes have the same size and therefore $d|n$ and also
 $n/d$ is even. This explains, for instance, why for a cycle $C$, $\chi_{i}(C) = 2$ if the length of $C$ is 0 mod 4, and
 $\chi_{i}(C) = 3$ otherwise. Properties of $d$-regular graphs which admit such colorings are discussed at length in Woldar~\cite{W} (see also Lazebnik and Woldar~\cite{LW}), which include a number of applications
 to graph theory problems and constructions from group theory.

 \bigskip

 In this paper, we consider a more coding theoretic problem, namely strong colorings of $q$-ary cube graphs or {\em Hamming graphs}, as they will provide some
 natural extremal examples for the coupon coloring number and injective coloring number. The {\em Hamming Graph}, $H_{n,q}$, has vertex set $[q]^n := \{1,2,\dots,q\}^n$ and two $q$-ary $n$-tuples are adjacent if and only if they differ in one coordinate -- that is, they have Hamming distance $1$.  The following theorem will be proved in Section \ref{cubes}.

\thm\label{hamming}
Let $k,n,q,r \in \mathbb N$ and $H = H_{n,q}$.
\begin{center}
\begin{tabular}{lp{5in}}
{\rm 1.} & If $q = 2$ and $n = 2^r$, then $\chi_c(H) = \chi_{i}(H) = n$. \\
{\rm 2.} & If $q$ is a prime power and $n = (q^k - 1)/(q - 1) + 1$, then $\chi_i(H) = q^k$. \\
{\rm 3.} & If $p, q$ are prime powers with $q^k \leq p$, and $n = (p^r - 1)/(p - 1) \cdot (q^k - 1)/(q - 1)$, then
$\chi_i(H) \leq p^r$.
\end{tabular}
\end{center}
\xthm

For $q > 2$, we have that $H_{n,q}$ is $(q - 1)n$-regular, and since $(q - 1)n$ does not divide $q^n$,
$\chi_c(H_{n,q}) < \chi_{i}(H_{n,q})$. In this sense, Theorem \ref{hamming}.2 is best possible, and there is
no analogue for $q > 2$ of Theorem \ref{hamming}.1. For $q = 2$ and $q = 3$, \"{O}sterg\aa rd~\cite{O} (see also Du, Kim and Pardalos~\cite{D}) showed that
$\chi_{i}(H_{n,q}) \sim (q - 1)n$, but this is not known for any $q > 3$. Theorem \ref{hamming}.3 gives supporting evidence for
the following conjecture:

\conj
For all $q > 2$, $\chi_{i}(H_{n,q}) \sim \chi_c(H_{n,q}) \sim (q - 1)n$ as $n \rightarrow \infty$.
\xconj

The different components of Theorem \ref{hamming} are proved in Section \ref{cubes}. We finally remark that $\chi_{i}(G)$
may be quadratically large relative to $\Delta(G)$: for instance if $G$ is the Erd\H{o}s-R\'{e}nyi polarity graph~\cite{ER}
of a projective plane of order $q$, then $\Delta(G) = q + 1$ whereas $\chi_{i}(G)  = |V(G)| = q^2 + q + 1$ since every pair of
vertices of $G$ is joined by a path of length two.

\section{Coupon coloring cubes}\label{cubes}

{\bf Proof of Theorem \ref{hamming}.1.} We prove $\chi_c(H_{n,2}) = \chi_{i}(H_{n,2}) = n$ when $n = 2^r$ for some $r \in \mathbb N$.
Since $H_{n,2}$ is $n$-regular, it is sufficient to find a coloring of $H_{n,2}$ with $n$ colors such that the neighborhood of each vertex contains all $n$ colors. To find such a coloring, identify $V(H_{n,2})$ with the power set of $\mathbb{F}_{2^r}$ in the natural way.  Specifically, order the elements of $\mathbb{F}_{2^r}$ in any way and for $v \in V(H_{n,2})$, identify $v$ with the subset of $\mathbb{F}_{2^r}$ that contains $x\in \mathbb{F}_{2^r}$ if and only if $v$ has a $1$ in the position associated with $x$. Now let $v$ be associated with subset $A$.  We also identify $2^r$ colors with the elements of $\mathbb{F}_{2^r}$, and we color vertex $v$ with color
\[
\sum_{x\in A} x.
\]
Now, if $v\sim w$, then $v$ and $w$ have Hamming distance $1$, which means that there is a unique $y\in \mathbb{F}_{2^r}$ such that $w$ is colored with color
\[
y+\sum_{x\in A} x.
\]
Note that we are using here that $\mathbb{F}_{2^r}$ has characteristic $2$. Then if $w_1$ and $w_2$ are distinct neighbors of $v$, the indices where they differ from $v$ must be distinct, and thus they must be colored with different colors. This proves Theorem \ref{hamming}.1. \qed

\bigskip

{\bf Proof of Theorem \ref{hamming}.2.} We prove that if $q$ is a prime power and
$n = (q^k - 1)/(q - 1)$, then $\chi_{i}(H_{n,q}) = q^k = (q - 1)n + 1$. To find such a injective coloring, we use a construction from coding theory~\cite{VW}. In $\mathbb{F}_q^k$, there are $n$ distinct $1$-dimensional subspaces.  Let $A$ be the $k\times n$ matrix with entries in $\mathbb{F}_q$ where each column is a representative from a distinct $1$-dimensional subspace.   Now let
\[
\mathcal{C} = \{x\in \mathbb{F}_q^n|Ax = 0\}.
\]
This is an $n-k$ dimensional subspace of $\mathbb{F}_q^n$.  Each $x\in \mathbb{F}_q^n$ represents a vertex of $H_{n,q}$. Let $H^2_{n,q}$ be the graph with vertex set $V(H_{n,q})$ where two vertices are adjacent if they are at the ends of a path of length $2$ in $H_{n,q}$ (in this case, this is the same as the square of $H_{n,q}$).  We claim that $\mathcal{C}$ forms an independent set in $H^2_{n,q}$.  In other words, $\mathcal{C}$ is a code with minimum distance $3$.  To see this, assume there are two vectors in $\mathcal{C}$ with distance less than $3$.  Without loss of generality, let one be the zero vector, and let the other be of the form $x^T=(a,b,0,...,0)$ where $a$ and $b$ are not both equal to $0$.  Then $Ax = 0$ implies that the first two columns are linearly dependent, a contradiction.

\medskip

Now, $\mathcal{C}$ is a subgroup of the additive group $\mathbb{F}_q^n$.  Since each of its cosets is also of minimum distance $3$, we can partition $\mathbb{F}_q^n$ into $q^k$ independent sets.  This gives that
\[
\chi(H^2_{n,q}) \leq q^k.
\]
(In fact, there is equality, since these independent sets are maximal by the sphere packing bound).  Since we have properly colored $H^2_{n,q}$, any pair of vertices in a neighborhood of a given vertex must have distinct colors.  This means that every vertex in $H_{n,q}$ sees $n(q-1) = q^k - 1$ distinct colors. This proves  Theorem \ref{hamming}.2.\qed

\bigskip

{\bf Proof of Theorem \ref{hamming}.3.} Let $n_q = (q^k - 1)/(q - 1)$ and $n_p = (p^r - 1)/(p - 1)$. We have to show that if $p,q$ are prime powers with $q^k \leq p$, and $m = n_p n_q$, then $\chi_{i}\left(H_{m,q}\right) \leq p^r$. Let $H_{n_p,p}$ and $H_{n_q,q}$ be strongly $p^r$-coupon colored and strongly $q^k$-coupon colored with color classes $\{P_1,P_2,\dots,P_{p^r}\}$ and $\{Q_1,Q_2,\dots,Q_{q^k}\}$ respectively, which we can do by Theorem \ref{hamming}.2.  We will create a strong $p^r$-coupon coloring of $H = H_{n_pn_q,q}$ with color classes $\{C_1,C_2,\dots,C_{p^r}\}$.  Let $(p_{j,1},...,p_{j,n_p})$ be a vertex in color class $P_j$.  In color class $C_j$ we place the following vertices of $V(H)$.
\[
\{(c_1\sqcup \cdots \sqcup c_{n_p}): c_i \in Q_{p_{j,i}}\}
\]
where $\sqcup$ represents concatenation. That is, we are replacing symbol $x\in [p]$ with all strings in color class $Q_x$.  Then because $Q_i$ and $P_j$ are independent sets for any $i$ and $j$, any pair of strings in $C_l$ has minimum distance $3$, which means $C_l$ is an independent set in $H^2$.  Since the $C_l$'s partition the $q$-ary strings of length $n_pn_q$, we have properly colored $H^2$, and thus have strongly $p^r$ colored $H$.
This proves Theorem \ref{hamming}.3. \qed

\section{Probabilistic tools}

The proof of Theorem \ref{main} is entirely probabilistic. To prove the lower bound $\chi_c(G) \geq (1 - \delta)d/\log d$ for $d$-regular graphs,
we use concentration inequalities and the local lemma:

\thm \label{local}
{\rm \cite{lll}} {\bf (Lov\' asz Local Lemma.)}
Let $A_i$ be a finite collection of events with $P(A_i) \le p$ for some $p$.
Suppose also that each $A_i$ is independent of all but at most $d$ other events.
Then if
\[ ep(d+1) \le 1,\]
with positive probability none of the $A_i$s occur.
\xthm

A sum of Bernoulli random variables is concentrated around its expectation via the following theorem.

\thm {\rm \cite{cher}} {\bf (Chernoff Bound.)}
Let $X_i$ be independent Bernoulli random variables and write $X := \sum X_i$.
Then for any $\gamma > 0$,
\[ \mathbb P( |X - E(X)| > \gamma E(x)) \le 2\exp\Bigl(- \f {\gamma^2 \mu}{2}\Bigr).\]
\xthm

When a martingale is Lipschitz, we can easily control the sum of the differences (using, say, the Hoeffding-Azuma Inequality).
If the martingale has a very low probability of failing to be Lipschitz, we can still exert some control over the sum of differences, though our tools need to be more precise.
This inequality appears in a slightly less general form in Shamir and Spencer~\cite{SS}, and in this form it was first stated and proved in the paper of Kim~\cite{K} using stopping times, and then restated in Chalker, Godbole, Hiczenko, Radcliffe and Ruehr~\cite{exceptionalprob} with a different proof.

\thm {\rm \cite{SS,K,exceptionalprob}}\label{exceptionalprob}
Let $(Z_i)_{i = 0}^m$ be a martingale such that $\triangle_i := Z_{i + 1} - Z_i$ satisfies $|\triangle_i| \leq B$ and let $c = (c_1,c_2,\dots,c_m) \in \mathbb R_+^m$ and $\|c\|^2 = \sum_{i = 1}^m c_i^2$. Then for $\lambda > 0$
and $\sigma \in (0,1)$,
\[ \mathbb P(|Z_m - Z_0| \geq \lambda) \leq e^{-\lambda^2 \sigma^2/8\|c\|^2} + (1 + \tfrac{B}{(1 - \sigma)\lambda}) \sum_{i = 1}^{m-1} \mathbb P(\triangle_i < -c_i).\]
\xthm

We will require the number of $d$-regular $n$-vertex graphs:

\thm {\rm \cite{mckay}} \label{thm:numreg}
The number of $d$-regular graphs on $n$ vertices with $nd$ even is
\[ (1+o(1))e^{-\f 14 (d^2-1)} \f {(dn)!}{{\f {dn}2}!2^{\f {dn}2} (d!)^{n}}.\]
\xthm

\section{Proof of Theorem \ref{main}}

{\bf Coupon coloring random regular graphs.} To prove the bound for random regular graphs in Theorem \ref{main}, let $G_{n,p}$ denote the Erd\H{o}s-R\'{e}nyi model of random graphs, where each
edge has probability $p = d/n$, independent of all other edges. Let $A_d$ be the event that $G \in G_{n,p}$ is $d$-regular, where
$n$ is even. By Theorem \ref{thm:numreg} and Stirling's Formula:
\begin{align*}
\mathbb P(A_d) &\sim  \frac{\sqrt{2}e^{\frac{1}{4}d^2 + \frac{1}{2}d - \frac{1}{4}}}{ (2 \pi d)^{n/2}}.
\end{align*}
Let $B$ be the event that $G \in G_{n,p}$ has $\chi_c(G) \geq k$
where $d = (1-\epsilon)k\log k$ -- we assume all quantities are integers, without affecting the forthcoming asymptotic computations.
For $1 \le i < j \le k$, let $E_{ij}$ be the event that every vertex in $V_i$ has a neighbor of color $j$. Then
\[ \mathbb P(B) \le \prod_{i < j} P(E_{ij})\]
since the $E_{ij}$ are independent events. Now the expression on the right is easily seen to be maximized when
all the color classes have size $n/k$ (we may assume $k|n$, since this does not affect the asymptotics of the
final answer). Then
\[ \log \mathbb P(E_{ij}) = \log (1 - (1 - p)^{n/k})^{n/k} \sim -nk^{\epsilon-2}.\]
We therefore have
\begin{align*}
\log \mathbb P(B) & \le \log \prod_{i < j} P(E_{ij}) \\
&\lesssim -nk^{\epsilon-2} \cdot {\binom k2} \\
&\lesssim -\frac{1}{2}nk^{\epsilon}.
\end{align*}
Since there are $k^n$ colorings,
\[ \log \mathbb P(\chi_c(G_{n,p}) \geq k) \lesssim n\log k - \frac{1}{2}nk^{\epsilon} \lesssim - \frac{1}{2}nk^{\epsilon}.\]
Finally, using the asymptotic formula for $\mathbb P(A_d)$,
\[ \log \mathbb P(B|A_d) \leq \log \mathbb P(B) - \log \mathbb P(A_d) \lesssim - \frac{1}{2}nk^{\epsilon} + \frac{1}{2}n\log d \lesssim -\frac{1}{2}nk^{\epsilon}\]
as required for Theorem \ref{main}. \qed

\bigskip

{\bf Coupon coloring regular graphs.} Here we prove the first statement in Theorem \ref{main}: we prove that for $\delta > 0$, as $d \rightarrow \infty$, every $d$-regular graph $G$ has
\[ \chi_{c}(G) \geq (1 - \delta)\frac{d}{\log d}.\]
The strategy involves two rounds of random coloring.
First, reserve a small subset $U \subset V$ of vertices, and randomly color $V \sm U$.
In the second round, identify the ``useful colors'' for each $u \in U$ and randomly color each $u$ from this smaller list.
For $u \in U$, the set of useful colors is precisely the set of colors $i$ such that there exists a neighbor $v$ of $u$
none of whose neighbors is colored $i$ -- in other words, $i$ is a missing color in the neighborhood of $v$.
The idea is that in the second round of coloring, since each vertex $v$ has many neighbors in $U$, it is very likely
that the colors missing in the neighborhood of $v$ in the first round now appear on neighbors of $v$ under the second round of coloring.

\medskip

{\bf The set $U$.} Let $d = (1 + \delta)k\log k$ and let $G = (V,E)$ be a $d$-regular graph.
Let $\eta < \dta$. Let $U$ be a randomly chosen subset of $V$ where the vertices are sampled independently with probability $d^{-\eta}$.
Let $U_v = U \cap \Gamma(v)$ for $v \in V$. Then $E(|U_v|) = d^{1 - \eta}$ for all $v \in V$.
Let $B_v$ be the event $||U_v| - d^{1-\eta}| > \gamma d^{1-\eta}$ where $\gamma = 4\sqrt{\log d/d^{1 - \eta}}$.
By the Chernoff Bound, for $v \in V$
\begin{align*}
\mathbb P(B_v) &< 2e^{-\gamma^2 d^{1 - \eta}/2} \\
&= 2e^{-8\log d} \\
&< d^{-4}.
\end{align*}
Now, the random variables $|U_v|, |U_u|$ are independent except when $v,u$ are adjacent or have a common neighbor in $V$; hence the dependency graph of
the events $B_v$ has degree at most $d^2$.
It follows from the local lemma, Theorem \ref{local}, that we can pick a specific $U \subset V$ such that
none of the events $B_v$ occurs, in other words,
\begin{equation}\label{u}
(1 - \gamma)d^{1 - \eta} < |U_v| < (1 + \gamma) d^{1 - \eta}
\end{equation}
holds simultaneously for every $v \in V$.

\bigskip

{\bf First round of coloring.} Randomly and independently assign a color from $[k]$ to each vertex of $V \sm U$. For $v \in V$, we
say that {\em $v$ misses the color $i$} if no vertex of $\Gamma(v) \sm U_v$ was assigned color $i$. Let $A_{v,i}$ denote this event.
Recall that $d = (1 + \delta)k\log k$ and so $(1 + \gamma)d^{1 - \eta} \leq k \log 2$ for large enough $d$. By (\ref{u}), for large enough $d$,
\begin{align*}
\mathbb P(A_{v,i}) &\le \left(1 - \tfrac{1}{k}\right)^{d - (1 + \gamma)d^{1 - \eta}} \\
&\leq \left(1 - \tfrac{1}{k}\right)^{d - k\log 2} \\
&\leq 2e^{-d/k} \; \; = \; \; 2k^{-1-\delta}.
\end{align*}
Now, for $u \in U$, consider the list $K_u$ of colors missed by at least one neighbor of $u$:
\[ K_u = \{i \in [k] : \exists v \in \Gamma(u), A_{v,i} \mbox{ occurs}\}.\]
For a particular color $i$, the union bound gives
\[ \mathbb P(i \in K_u) \le \sum_{v \sim u} \mathbb P(A_{v,i}) \le 2dk^{-1-\dta}\]
and hence
\[ \mathbb E(|K_u|) = \sum_{i \in [k]} \mathbb P(i \in K_u) \leq 2dk^{-\delta}.\]
Fix $u \in U$ and let $\Gam(u) = \{ v_1,\ldots, v_d\}$. Let
\[ Z_0 = \mathbb E(|K_u|), \quad Z_i = \mathbb E \big( |K_u| \, \big | \, \Gam(v_1),\ldots, \Gam(v_i)\big ).\]
That is, $Z_i$ is an exposure martingale on $|K_u|$ where at step $i$ we reveal the colors of the neighbors of $v_i$.
Note that if $|Z_i - Z_{i-1}| > c$, then $v_i$ necessarily misses more than $c$ colors.
This happens with probability at most
\begin{align*}
\binom kc \cdot 2k^{-c(1 + \delta)} &\le 2 k^ck^{-c(1+\dta)} \\
&\le 2k^{-\dta c}.
\end{align*}
Therefore $(Z_j)_{j = 0}^d$ is $c$-Lipschitz with exceptional probability at most $2k^{-\dta c}$.

We now use Theorem \ref{exceptionalprob} to obtain concentration of $|K_u|$ around its expectation.
Note that the martingale difference $|Z_j-Z_{j-1}|$ is obviously bounded by $d$.
Choosing $\lam = 8c\sqrt {d\log d}$ and $\sig = 9/10$ we obtain
\begin{align*}
\mathbb P({|K_u| > 2dk^{-\dta} + \lambda})
&\le \mathbb P(||K_u| - \mathbb E(|K_u|)| > \lambda) \\
&\le 2 \exp \Bigl({-8\sigma^2 \log d}\Bigr) + d \Bigl(1 + \f {10d}{\lambda}\Bigr)2k^{-\dta c} \\
&\le 2 d^{-6} + d^2 k^{-\dta c}.
\end{align*}
Now setting $c = 8\dta \inv$ shows that for large enough $d$,
\begin{align*}
\mathbb P(|K_u| > 2dk^{-\dta} + 64\dta\inv\sqrt{d\log d})
&\le 2 d^{-6} + d^2 k^{-8} \\
&\le  2d^{-6} + d^{-5} \\
&\le d^{-4}
\end{align*}
The random variables $|K_u|, |K_v|$ are dependent if and only if $u,v \in U$ share a common neighbor in $V$ or are adjacent.
Hence the dependency graph has degree at most $d^2$, and an application of the Local Lemma shows that there is a coloring for which
\[ |K_u| < 2dk^{-\dta} + 64\dta\inv\sqrt{d\log d} \leq 4dk^{-\delta}\]
for every $u \in U$ and large enough $d$.

\bigskip

{\bf Second round of coloring.} Finally, assign each $u \in U$ a color chosen uniformly at random from $K_u$.
By (\ref{u}), the probability that $v\in V$ does not have color $i$ in $\Gamma(v)$ is bounded by
\begin{align*}
\prod_{u \in U_v} \Bigl(1-\f {1}{|K_u|}\Bigr) &\le \Bigl({1 - \f 1{4dk^{-\dta}}}\Bigr)^{(1-\gamma)d^{1-\eta}} \\
&\le \exp\Bigl( - \f {(1-\gamma)d^{1-\eta}}{4dk^{-\dta}}\Bigr) \\
&= \exp\Bigl(- \f {k^\dta (1-\gamma)}{4d^\eta}\Bigr).
\end{align*}
The union bound shows that the event that $v$ does not have all $k$ colors in its neighborhood is at most
\begin{align*}
k\exp\Bigl(-\f {k^\dta (1-\gamma)}{4d^\eta}\Bigr) &= \exp\Bigl(-\f {k^{\dta-\eta} (\log k)^{1-\eta} (1-\gamma)}{4(1+\dta)^\eta}\Bigr) \\
&\le d^{-4},
\end{align*}
since $\dta-\eta > 0$.
Again, these events are independent except when two vertices are adjacent or share a common neighbor, so the dependency graph has degree at most $d^2$.
A final application of the local lemma, Theorem \ref{local}, shows that with positive probability, none of the events $A_{v,i}$ occur
for $v \in V$ and $i \in [k]$. \qed

\section{Concluding Remarks}

$\bullet$ We showed that every $d$-regular graph $G$ has $\chi_c(G) \geq (1 - o(1))d/\log d$, and random regular graphs demonstrate the
tightness of this result. It would be interesting to give some explicit examples for each $\delta > 0$ and each $d$ of infinitely many $d$-regular graphs for which $\chi_c(G) \leq (1 + \delta)d/\log d$. In particular, it seems likely that $\chi_c(G_q) \leq (1 + o(1))d/\log d$ when $G_q$ is the Paley graph, and $d = (q - 1)/2$. One may also ask what is $\chi_c(G)$ when $G$ is the bipartite incidence graph of a projective plane of order $q$.

\bigskip

$\bullet$ We leave the open question of showing $\chi_{i}(H_{n,q}) \sim q(n - 1)$ and $\chi_c(H_{n,q}) \sim q(n - 1)$ for $q > 2$.
More generally, if $G$ is a graph it would be interesting to determine any relations between $\chi_c(G)$ and $\chi_c(G \times G)$,
where $G \times G$ denotes the Cartesian product of $G$ with itself.

\end{document}